\documentclass[titlepage,12pt]{article}

\usepackage{graphicx}

\usepackage{amssymb}

\newcommand{\R}{\mathbb R}

\begin{document}

\baselineskip=18pt

\begin{center}
{\Large{\bf Bifurcation Analysis of the \\ Watt Governor System}}
\end{center}

\begin{center}
{\large Jorge Sotomayor}
\end{center}
\begin{center}
{\em Instituto de Matem\'atica e Estat\'{\i}stica, Universidade de
S\~ao Paulo\\ Rua do Mat\~ao 1010, Cidade Universit\'aria\\ CEP
05.508-090, S\~ao Paulo, SP, Brazil
\\}e--mail:sotp@ime.usp.br
\end{center}

\begin{center}
{\large Luis Fernando Mello}
\end{center}
\begin{center}
{\em Instituto de Ci\^encias Exatas, Universidade Federal de
Itajub\'a\\Avenida BPS 1303, Pinheirinho, CEP 37.500-903,
Itajub\'a, MG, Brazil
\\}e--mail:lfmelo@unifei.edu.br
\end{center}

\begin{center}
{\large Denis de Carvalho Braga}
\end{center}
\begin{center}
{\em Instituto de Sistemas El\'etricos e Energia, Universidade
Federal de Itajub\'a\\Avenida BPS 1303, Pinheirinho, CEP
37.500-903, Itajub\'a, MG, Brazil
\\}e--mail:braga\_denis@yahoo.com.br
\end{center}

\begin{center}
{\bf Abstract}
\end{center}

\vspace{0.1cm}

This paper pursues the study carried out by the authors in {\it
Stability and Hopf bifurcation in the Watt governor system}
\cite{smb}, focusing on the codimension one Hopf bifurcations in
the centrifugal Watt governor differential system, as presented in
Pontryagin's book {\it Ordinary Differential Equations},
\cite{pon}. Here are studied the codimension two and three Hopf
bifurcations and the pertinent Lyapunov stability coefficients and
bifurcation diagrams, illustrating the number, types and positions
of bifurcating  small amplitude periodic orbits, are determined. As
a consequence it is found a region in the space of parameters
where an attracting  periodic orbit coexists with an attracting
equilibrium.

\vspace{0.1cm}



\newpage
\baselineskip=20pt

\section{\bf Introduction}\label{intro}

The Watt centrifugal governor is a device that automatically
controls the speed of an engine. Dating to 1788, it can be taken
as the starting point for the theory  of automatic control (see
MacFarlane \cite{mac} and references therein). In this paper the
system coupling the Watt-centrifugal-governor and the steam-engine
will be called simply the Watt Governor System (WGS). See Section
\ref{watt} for a description and illustration, in Fig.
\ref{wattgov}, of this system.

Landmarks for the study of the local stability analysis of the WGS
are the works of Maxwell \cite{max} and Vyshnegradskii
\cite{vysh}. A simplified version of the WGS local stability based
on the work of Vyshnegradskii is presented by Pontryagin
\cite{pon}. A local stability study  generalized to a more general
Watt governor design was carried out by  Denny \cite{denny}  and
pursued by the authors in \cite{smb}.

Enlightening historical comments about the Watt governor local
mathematical stability and oscillatory analysis can be found in
MacFarlane \cite{mac} and Denny \cite{denny}. There, as well as in
\cite{pon}, we learn that toward the mid $XIX$ century,
improvements in the engineering design led to  less reliable
operations in the WGS, leading to fluctuations and oscillations
instead of the ideal  stable constant speed  output  requirement.
The first mathematical analysis of the stability conditions and
subsequent indication of the modification in the design to avoid
the problem was carried out by  Maxwell \cite{max} and, in a user
friendly style likely to be better understood by engineers, by
Vyshnegradskii \cite{vysh}.

From the mathematical point of view, the oscillatory, small
amplitude, behavior in the WGS can be associated to a periodic
orbit that appears from a Hopf bifurcation. This was established
by Hassard et al. in \cite{has1} and  Al-Humadi and Kazarinoff in
\cite{humadi}. Another procedure, based in the method of harmonic
balance,  has been suggested by Denny \cite{denny}  to detect
large amplitude oscillations.

In \cite{smb} we characterized the surface of Hopf bifurcations in
a WGS, which is more general than that presented by Pontryagin
\cite{pon}, Al-Humadi and Kazarinoff \cite{humadi} and Denny
\cite{denny}. See Theorem \ref{teohopf} and Fig. \ref{curvag=0}
for a review of the critical curve on the  surface where the first
Lyapunov coefficient vanishes.

In the present paper, restricting ourselves to  Pontryagin's
system, we go deeper investigating the stability of the
equilibrium along the above mentioned critical curve. To this end
the second Lyapunov coefficient is calculated (Theorem
\ref{teohopf2}) and it is established that it vanishes at a unique
point (see Fig. \ref{curvagh} and \ref{pointsTR}).  The third
Lyapunov coefficient is calculated at this point (Theorem
\ref{teohopf3}) and found to be positive. The pertinent
bifurcation diagrams are established. See Fig. \ref{pointT},
\ref{pointR} and \ref{pointR1}. A conclusion  derived from these
diagrams, concerning the region ---a solid ``tongue"--- in the
space of parameters where an attracting periodic orbit coexists
with an attracting equilibrium, is specifically commented in
Section \ref{conclusion}.

The extensive calculations involved in Theorems \ref{teohopf2} and
\ref{teohopf3} have been corroborated with the software
MATHEMATICA 5 \cite{math} and the main steps have been posted in
the site \cite{mello}.

This paper is organized as follows. In Section \ref{watt} we
introduce the WGS and review the Pontryagin differential equations
\cite{pon}. The stability of the equilibrium points is also
analyzed. This section is essentially a review of \cite{pon, has1,
humadi, smb}.  The Hopf bifurcations in the WGS differential
equations are studied in Sections \ref{codim} and \ref{hopf}.
Expressions for the second and third Lyapunov coefficients, which
fully clarify their sign, are obtained, pushing forward the method
found in the works of Kuznetsov \cite{kuznet, kuznet2}. With this
data, the bifurcation diagrams are established. Concluding
comments, synthesizing and interpreting the results achieved here,
are presented in Section \ref{conclusion}.

\section{The  Watt centrifugal governor system}\label{watt}

\newtheorem{teo}{Theorem}[section]
\newtheorem{lema}[teo]{Lemma}
\newtheorem{prop}[teo]{Proposition}
\newtheorem{cor}[teo]{Corollary}
\newtheorem{remark}[teo]{Remark}

\subsection{Differential equations for the Watt governor system }\label{diffequat}

\begin{figure}[!h]

\centerline{
\includegraphics[width=13cm]{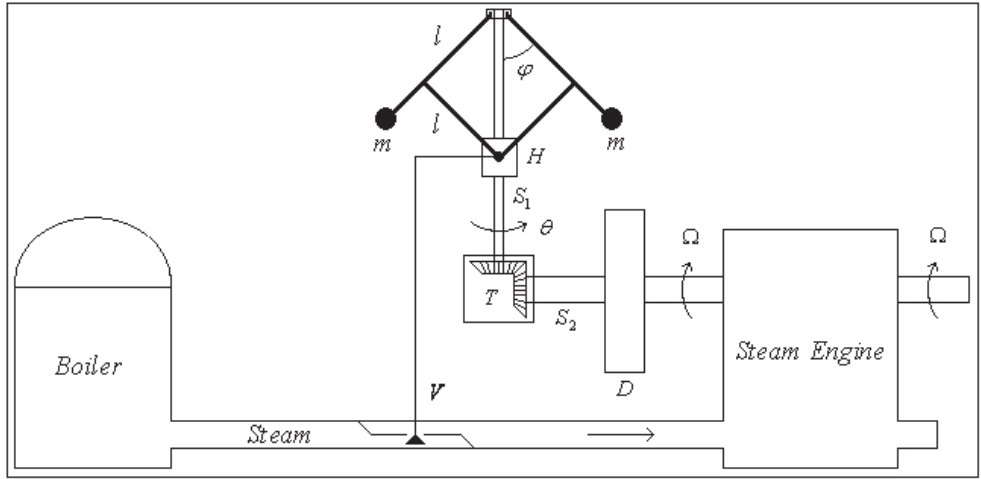}}

\caption{{\small Watt-centrifugal-governor-steam-engine system}.}

\label{wattgov}

\end{figure}

According to Pontryagin \cite{pon}, p. 217, the differential
equations of the WGS illustrated in Fig. \ref{wattgov} are
\begin{eqnarray}\label{wattde}
\frac{d \; \varphi}{d \tau} &=& \psi \nonumber\\
\frac{d \; \psi}{d \tau} &=& c^2 \; \Omega^2 \; \sin \varphi \;
\cos \varphi -
\frac{g}{l} \; \sin \varphi - \frac{b}{m} \; \psi \\
\frac{d \; \Omega}{d \tau} &=& \frac{1}{I} \; \left( \mu \cos
\varphi -F \right) \nonumber
\end{eqnarray}
\noindent where $\varphi \in \left( 0,\frac{\pi}{2} \right)$ is
the angle of deviation of the arms of the centrifugal governor
from its vertical axis $S_1$, $\Omega \in [0,\infty)$ is the
angular velocity of the rotation of the flywheel $D$, $\theta$ is
the angular velocity of $S_1$, $l$ is the length of the arms, $m$
is the mass of each ball, $H$ is a  sleeve which supports the arms
and slides along $S_1$, $T$ is a set of transmission gears, $V$ is
the valve that determines the supply of steam to the engine,
$\tau$ is the time, $\psi = d \varphi/d \tau$, $g$ is the standard
acceleration of gravity, $\theta = c \: \Omega $, $c > 0$ is a
constant transmission ratio, $b > 0$ is a constant of the
frictional force of the system, $I$ is the moment of inertia of
the flywheel, $F$ is an equivalent torque of the load and $\mu >
0$ is a proportionality constant. The reader is referred  to
Pontryagin  \cite{pon} for  the derivation  of (\ref{wattde}) from
Newton's Second Law of Motion.

After the following change in the coordinates and time
\begin{equation}
x = \varphi, \:\: y = \sqrt {\frac {l}{g}} \: \psi, \:\: z = c \:
\sqrt {\frac {l}{g}} \: \Omega, \:\: \tau = \sqrt {\frac {l}{g}}
\: t, \label{mudanca}
\end{equation}
\noindent the differential equations (\ref{wattde}) can be written
as
\begin{eqnarray}\label{wattdem}
x' = \frac{d x}{d t} &=& y \nonumber\\
y' = \frac{d y}{d t} &=& z^2 \; \sin x \; \cos x -
\sin x - \varepsilon \; y \\
z' = \frac{d z}{d t} &=& \alpha \; \left( \cos x - \beta \right)
\nonumber
\end{eqnarray}
\noindent where $\alpha > 0$, $0 < \beta < 1$ and $\varepsilon
> 0 $, given by
\begin{equation}
\varepsilon = \frac{b}{m} \; \sqrt \frac{l}{g}, \:\: \alpha =
\frac{c \; l \; \mu}{g \; I}, \:\: \beta = \frac{F}{\mu},
\label{parameters}
\end{equation}
are the normalized variable parameters. Thus the differential
equations (\ref{wattdem}) are in fact a three-parameter family of
differential equations which can be rewritten as ${\bf x}' = f
({\bf x}, {\bf \mu})$, where
\begin{equation}
{\bf x} =(x,y,z) \in \left( 0, \frac{\pi}{2} \right) \times \R
\times [0,\infty), \: \: \: {\bf \mu} = (\beta, \alpha,
\varepsilon) \in \left( 0, 1\right) \times \left( 0, \infty
\right) \times \left(0,\infty \right)
 \label{paramfase}
\end{equation}
and
\begin{equation}
f({\bf x},{\bf \mu}) = \left( y, z^2 \; \sin x \; \cos x - \sin x
- \varepsilon \; y, \alpha \; \left( \cos x - \beta \right)
\right). \label{wattdemo}
\end{equation}

\subsection{Stability analysis of the equilibrium
points}\label{stability}

The differential equations (\ref{wattdem}) have one admissible
equilibrium point
\begin{equation}
P_0 = (x_0,y_0,z_0) =  \left( \arccos \beta, 0, \sqrt
\frac{1}{\beta} \right). \label{equilibrium}
\end{equation}
The Jacobian matrix of $f$ at $P_0$ has the form
\begin{equation}
Df \left( P_0 \right) = \left( \begin{array}{ccc}
0 & 1 & 0 \\
\\-\displaystyle \frac{1- \beta^2}{\beta} & - \varepsilon & 2 \sqrt{\beta (1-\beta^2)} \\
\\-\displaystyle \alpha \sqrt{1 - \beta^2}    & 0  & 0
\end{array} \right)
\label{jacobian}
\end{equation}
\noindent and its characteristic polynomial is given by
$p(\lambda)$, with
\begin{equation}
-p(\lambda)= \lambda^3 + \varepsilon \: \lambda^2 + \frac{1-
\beta^2}{\beta} \: \lambda + 2\: \alpha \: {\beta}^{3/2} \:
\frac{1-\beta^2}{\beta}. \label{charac}
\end{equation}

\begin{figure}[!h]
\centerline{
\includegraphics[width=8cm]{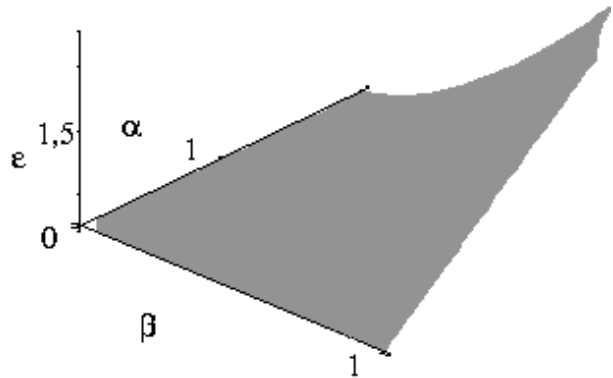}}

\caption{{\small Surface of critical parameters $\varepsilon_c = 2
\; \alpha \; {\beta}^{3/2}$}.}

\label{criticalsurf}
\end{figure}

\begin{teo}
For all
\begin{equation}
\varepsilon > 2 \: \alpha \: \beta ^{3/2} \label{asymp}
\end{equation}
the WGS differential equations (\ref{wattdem}) have an
asymptotically stable equilibrium point at $P_0$. If $0 <
\varepsilon < 2 \; \alpha \; \beta^{3/2}$ then $P_0$ is unstable.
\label{teorema1}
\end{teo}

The proof of this theorem can essentially be found in Pontryagin
\cite{pon}; it has also been established in a more general setting
in \cite{smb}.

The surface of critical parameters $\mu_0 = (\beta,
\alpha,\varepsilon_c)$ such that $\varepsilon_c = \varepsilon
(\beta,\alpha) = 2 \: \alpha \: \beta ^{3/2}$ is illustrated in
Fig. \ref{criticalsurf}. In the Section \ref{hopf} we will analyze
the stability of $P_0$ as $\varepsilon_c = 2 \: \alpha \: \beta
^{3/2}$. The change in the stability at the equilibrium $P_0$ as
the parameters cross the critical surface produces a Hopf
bifurcation in the WGS, whose analysis has been carried out by
\cite{humadi}, \cite{has1} and, in a more general setting, by
\cite{smb}.

From  (\ref{parameters}), $\varepsilon$ represents the friction
coefficient of the system. The case $\varepsilon =0$ maybe of
theoretical interest  due to  its connection  with conservative
systems.  However, as made explicit in {\it Vyshnegradskii's
Rules},  friction is an essential ingredient to attain stability.
This point is neatly presented in Pontryagin \cite{pon}, of which
Figure \ref{criticalsurf} is a geometric, dimensionless,
synthesis.

\section{Lyapunov coefficients}\label{codim}

The beginning of this section is a review of the method found in
\cite{kuznet}, pp 177-181, and  in  \cite{kuznet2} for the
calculation of the first and second Lyapunov coefficients. The
calculation of the third Lyapunov coefficient has not been found
by the authors in the current literature. The extensive
calculations and the  long expressions for these coefficients have
been corroborated with the software MATHEMATICA 5 \cite{math}.

Consider the differential equations
\begin{equation}
{\bf x}' = f ({\bf x}, {\bf \mu}), \label{diffequat}
\end{equation}
where ${\bf x} \in \R^n$ and ${\bf \mu} \in \R^m$ are respectively
vectors representing  phase variables and control parameters.
Assume that $f$ is of class $C^{\infty}$ in $\R^n \times \R^m$.
Suppose (\ref{diffequat}) has an equilibrium point ${\bf x} = {\bf
x_0}$ at ${\bf \mu} = {\bf \mu_0}$ and, denoting the variable
${\bf x}-{\bf x_0}$ also by ${\bf x}$, write
\begin{equation}
F({\bf x}) = f ({\bf x}, {\bf \mu_0}) \label{Fhomo}
\end{equation}
as
\begin{eqnarray}\label{taylorexp}
F({\bf x}) = A{\bf x} + \frac{1}{2} \: B({\bf x},{\bf x}) +
\frac{1}{6} \: C({\bf x}, {\bf x}, {\bf x}) + \: \frac{1}{24} \:
D({\bf x}, {\bf x}, {\bf x}, {\bf x}) {\nonumber} \\+
\frac{1}{120} \: E({\bf x}, {\bf x}, {\bf x}, {\bf x}, {\bf x}) +
\frac{1}{720} \: K({\bf x}, {\bf x}, {\bf x}, {\bf x}, {\bf x},
{\bf x})\\+ \frac{1}{5040} \: L({\bf x}, {\bf x}, {\bf x}, {\bf
x}, {\bf x}, {\bf x}, {\bf x}) + O(|| {\bf x} ||^8){\nonumber},
\end{eqnarray}
\noindent where $A = f_{\bf x}(0,{\bf \mu_0})$ and
\begin{equation}
B_i ({\bf x},{\bf y}) = \sum_{j,k=1}^n \frac{\partial ^2
F_i(\xi)}{\partial \xi_j \: \partial \xi_k} \bigg|_{\xi=0} x_j \;
y_k, \label{Bap}
\end{equation}
\begin{equation}
C_i ({\bf x},{\bf y},{\bf z}) = \sum_{j,k,l=1}^n \frac{\partial ^3
F_i(\xi)}{\partial \xi_j \: \partial \xi_k \: \partial \xi_l}
\bigg|_{\xi=0} x_j \; y_k \: z_l, \label{Cap}
\end{equation}
\begin{equation}
D_i ({\bf x},{\bf y},{\bf z},{\bf u}) = \sum_{j,k,l,r=1}^n
\frac{\partial ^4 F_i(\xi)}{\partial \xi_j \: \partial \xi_k \:
\partial \xi_l \: \partial \xi_r} \bigg|_{\xi=0} x_j \;
y_k \: z_l \: u_r, \label{Dap}
\end{equation}
\begin{equation}
E_i ({\bf x},{\bf y},{\bf z},{\bf u},{\bf v}) = \sum_{j,k,l,r,p
=1}^n \frac{\partial ^5 F_i(\xi)}{\partial \xi_j \: \partial \xi_k
\: \partial \xi_l \: \partial \xi_r \: \partial \xi_p}
\bigg|_{\xi=0} x_j \; y_k \: z_l \: u_r \: v_p, \label{Eap}
\end{equation}
\begin{equation}
K_i ({\bf x},{\bf y},{\bf z},{\bf u},{\bf v},{\bf s}) =
\sum_{j,\ldots,q =1}^n \frac{\partial ^6 F_i(\xi)}{\partial \xi_j
\: \partial \xi_k \: \partial \xi_l \: \partial \xi_r \: \partial
\xi_p \: \partial \xi_q} \bigg|_{\xi=0} x_j \; y_k \: z_l \: u_r
\: v_p \: s_q, \label{Kap}
\end{equation}
\begin{equation}
L_i ({\bf x},{\bf y},{\bf z},{\bf u},{\bf v},{\bf s},{\bf t}) =
\sum_{j,\ldots,h =1}^n \frac{\partial ^7 F_i(\xi)}{\partial \xi_j
\partial \xi_k
\partial \xi_l  \partial \xi_r  \partial \xi_p  \partial \xi_q \partial \xi_h} \bigg|_{\xi=0} x_j \;
y_k \: z_l \: u_r \: v_p \: s_q \: t_h, \label{Lap}
\end{equation}
\noindent for $i = 1, \ldots, n$.

Suppose $({\bf x_0}, {\bf \mu_0})$ is an equilibrium point of
(\ref{diffequat}) where the Jacobian matrix $A$ has a pair of
purely imaginary eigenvalues $\lambda_{2,3} = \pm i \omega_0$,
$\omega_0 > 0$, and admits no other eigenvalue with zero real
part. Let $T^c$ be the generalized eigenspace of $A$ corresponding
to $\lambda_{2,3}$. By this is meant that it is the largest
subspace invariant by $A$ on which the eigenvalues are
$\lambda_{2,3}$.

Let $p, q \in \mathbb C ^n$ be vectors such that
\begin{equation}
A q = i \omega_0 \: q,\:\: A^{\top} p = -i \omega_0 \: p, \:\:
\langle p,q \rangle = \sum_{i=1}^n \bar{p}_i \: q_i \:\: = 1,
\label{normalization}
\end{equation}
where $A^{\top}$ is the transposed matrix. Any vector $y \in T^c$
can be represented as $y = w q + \bar w \bar q$, where $w =
\langle p , y \rangle \in \mathbb C$. The two dimensional center
manifold can be parametrized by $w , \bar w$, by means of an
immersion of the form  ${\bf x} = H (w, \bar w)$, where $H:\mathbb
C^2 \to \R^n$ has a Taylor expansion of the form
\begin{equation}
H(w,{\bar w}) = w q + {\bar w}{\bar q} + \sum_{2 \leq j+k \leq 7}
\frac{1}{j!k!} \: h_{jk}w^j{\bar w}^k + O(|w|^8), \label{defH}
\end{equation}
with $h_{jk} \in \mathbb C ^n$ and  $h_{jk}={\bar h}_{kj}$.
Substituting this expression into (\ref{diffequat}) we obtain the
following differential equation
\begin{equation} \label{ku}
H_w w' + H_{\bar w} {\bar w}' = F (H(w,{\bar w})),
\end{equation}
where $F$ is given by (\ref{Fhomo}).

The complex vectors $h_{ij}$ are to be determined so that system
(\ref{ku}), on the chart $w$ for a central manifold, writes as
follows
\[
w'= i \omega_0 w + \frac{1}{2} \; G_{21} w |w|^2 + \frac{1}{12} \;
G_{32} w |w|^4 + \frac{1}{144} \; G_{43} w |w|^6 + O(|w|^8),
\]
with $G_{jk} \in \mathbb C$.

Solving for the vectors $h_{ij}$ the  system of linear equations
defined by the coefficients of the quadratic terms of (\ref{ku}),
taking into account the coefficients of $F$ in the expressions
(\ref{taylorexp}) and (\ref{Bap}), one has
\begin{equation}
h_{11}=-A^{-1}B(q,{\bar q}) \label{h11},
\end{equation}
\begin{equation}
h_{20}=(2i\omega_0 I_n - A)^{-1}B(q,q),\label{h20}
\end{equation}
where $I_n$ is the unit $n \times n$ matrix. Pursuing the
calculation to cubic terms, from the coefficients of the terms
$w^3$ in (\ref{ku}) follows that
\begin{equation}
h_{30}=(3 i \omega_0 I_n -A)^{-1} \left[ 3B(q,h_{20})+C(q,q,q)
\right] \label{h30}.
\end{equation}

From the coefficients of the terms $w^2 {\bar w}$ in (\ref{ku})
one obtains a singular system for $h_{21}$
\begin{equation}
(i \omega_0 I_n -A)h_{21}=C(q,q,{\bar q})+B({\bar q},h_{20})+ 2
B(q,h_{11})-G_{21}q, \label{h21m}
\end{equation}
which has a solution if and only if
\[
\langle p, C(q,q,\bar q) + B(\bar q, h_{20}) + 2 B(q,h_{11})
-G_{21} q \rangle = 0.
\]
Therefore
\begin{equation}
G_{21}= \langle p, C(q,q,\bar q) + B(\bar q, (2i \omega_0 I_n
-A)^{-1} B(q,q)) - 2 B(q,A^{-1} B(q,\bar q)) \rangle.
\end{equation}

The {\it first Lyapunov coefficient} $l_1$ is defined by
\begin{equation}
l_1 =  \frac{1}{2} \: {\rm Re} \; G_{21}. \label{defcoef}
\end{equation}

The complex vector $h_{21}$ can be found by solving the
nonsingular $(n+1)$-dimensional system
\begin{equation}
\left( \begin{array}{cc}
i \omega_0 I_n -A & q \\
\\ {\bar p} & 0
\end{array} \right) \left( \begin{array} {c}
h_{21}\\
\\s \end{array} \right)= \left( \begin{array}{c} C(q,q,\bar q) +
B(\bar q, h_{20}) + 2 B(q,h_{11}) -G_{21} q \\
\\0 \end{array} \right),\label{h20}
\end{equation}
with the condition $\langle p, h_{21} \rangle = 0$.

For the sake of completeness, in Remark \ref{nonsingular} we prove
that the system (\ref{h20}) is nonsingular and that if $(v,s)$ is
a solution of (\ref{h20}) with the condition $\langle p, v \rangle
= 0$ then $v$ is a solution of (\ref{h21m}).

\begin{remark}
Write $\R^n = T^c \oplus T^{su}$, where $T^c$ and $T^{su}$ are
invariant by $A$. It can be proved that $y \in T^{su}$ if and only
if $\langle p , y \rangle = 0$. Define
\[
a = C(q,q,\bar q) + B(\bar q, h_{20}) + 2 B(q,h_{11}) -G_{21} q.
\]
Let $(v,s)$ be a solution of the homogeneous equation obtained
from (\ref{h20}). Equivalently
\begin{equation}
(i \omega_0 I_n -A)v + s q = 0, \: \: \langle p , v \rangle = 0.
\label{sistema1}
\end{equation}
From the second equation of (\ref{sistema1}), it follows that $v
\in T^{su}$, and thus $(i \omega_0 I_n -A)v \in T^{su}$. Therefore
$\langle p, (i \omega_0 I_n -A)v \rangle = 0$. Taking the inner
product of $p$ with the first equation of (\ref{sistema1}) one has
$\langle p, (i \omega_0 I_n -A)v + s q \rangle = 0$, which can be
written as $ \langle p, (i \omega_0 I_n -A)v \rangle + s \langle
p, q \rangle = 0$. Since $\langle p, q \rangle = 1$ and $\langle
p, (i \omega_0 I_n -A)v \rangle = 0$ it follows that $s = 0$.
Substituting $s = 0 $ into the first equation of (\ref{sistema1})
one has $(i \omega_0 I_n -A)v = 0 $. This implies that
\begin{equation}
v = \alpha q, \: \alpha \in \mathbb C. \label{alpha}
\end{equation}
But $ 0 = \langle p , v \rangle = \langle p , \alpha q \rangle =
\alpha \langle p, q \rangle = \alpha$. Substituting $\alpha =0$
into (\ref{alpha}) it follows that $v=0$. Therefore $(v,s) =
(0,0)$.

Let $(v,s)$ be a solution of (\ref{h20}). Equivalently
\begin{equation}
(i \omega_0 I_n -A)v + s q = a, \: \langle p , v \rangle = 0.
\label{sistema}
\end{equation}
From the second equation of (\ref{sistema}), it follows that $v
\in T^{su}$ and thus $(i \omega_0 I_n -A)v \in T^{su}$. Therefore
$ \langle p , (i \omega_0 I_n -A)v \rangle = 0$. Taking the inner
product of $p$ with the first equation of (\ref{sistema}) one has
$\langle p, (i \omega_0 I_n -A)v + s q \rangle = \langle p, a
\rangle$, which can be written as
\[
\langle p, (i \omega_0 I_n -A)v \rangle + s \langle p, q \rangle =
\langle p, a \rangle.
\]
As $\langle p, a \rangle = 0$, $\langle p, q \rangle = 1$ and
$\langle p, (i \omega_0 I_n -A)v \rangle = 0$ it follows that $s =
0$. Substituting $s = 0$ into the first equation of
(\ref{sistema}) results $(i \omega_0 I_n -A)v = a$. Therefore $v$
is a solution of (\ref{h21m}).

\label{nonsingular}
\end{remark}

The procedure above will be adapted below in connection with the
determination of $h_{32}$ and $h_{43}$.

From the coefficients of the terms $w^4$, $w^3 {\bar w}$ and $w^2
{\bar w ^2}$ in (\ref{ku}), one has respectively
\begin{equation}
h_{40}= (4i \omega_0 I_n -A)^{-1}[3B(h_{20},h_{20}) + 4
B(q,h_{30})+ 6 C(q,q,h_{20})+ D(q,q,q,q)],\label{h40}
\end{equation}
\begin{eqnarray}\label{h31}
h_{31}=(2i \omega_0 I_n -A)^{-1}[3 B(q,h_{21}) + B(\bar q, h_{30})
+ 3 B(h_{20},h_{11})  \nonumber \\ + 3 C(q,q,h_{11}) +3 C(q,\bar
q,h_{20}) + D(q,q,q,\bar q) - 3G_{21}h_{20}],
\end{eqnarray}
\begin{eqnarray}\label{h22}
h_{22}= -A^{-1} [D(q,q,\bar q, \bar q) + 4 C(q, \bar q, h_{11}) +
C(\bar q, \bar q, h_{20}) + C(q,q,{\bar h}_{20}) \nonumber \\
+ 2 B(h_{11},h_{11}) + 2 B(q,{\bar h}_{21}) + 2 B(\bar q, h_{21})+
B({\bar h}_{20},h_{20})],
\end{eqnarray}
where the term $-2h_{11}(G_{21}+{\bar G}_{21})$ has been omitted
in the last equation, since $G_{21}+{\bar G}_{21} = 0$ as $l_1 =
0$.

Defining $\mathcal H_{32}$ as
\begin{eqnarray}\label{H32}
\mathcal H_{32} = 6 B(h_{11},h_{21})+ B({\bar h}_{20},h_{30}) + 3
B({\bar h}_{21},h_{20})+ 3 B(q,h_{22}) \nonumber \\ + 2 B(\bar q,
h_{31}) +6 C(q,h_{11},h_{11}) + 3 C(q, {\bar h}_{20}, h_{20})+ 3
C(q,q,{\bar h}_{21}) \nonumber \\ +6 C(q,\bar q, h_{21}) + 6
C(\bar q, h_{20}, h_{11}) + C(\bar q, \bar q, h_{30}) +
D(q,q,q,{\bar h}_{20}) \\ + 6 D(q,q,\bar q,h_{11}) + 3
D(q, \bar q,\bar q, h_{20}) + E(q,q,q,\bar q,\bar q) \nonumber \\
-6 G_{21}h_{21} - 3 {\bar G}_{21} h_{21}, \nonumber
\end{eqnarray}
and from the coefficients of the terms $w^3 {\bar w}^2$ in
(\ref{ku}), one has a singular system for $h_{32}$
\begin{equation}
(i \omega_0 I_n -A)h_{32}= \mathcal H_{32} - G_{32}q,
\label{h32m}
\end{equation}
which has solution if and only if
\begin{equation}
\langle p, \mathcal H_{32} - G_{32}q \rangle = 0, \label{H32m}
\end{equation}
where the terms $-6 G_{21}h_{21} - 3 {\bar G}_{21} h_{21}$ in the
last line of (\ref{H32}) actually does not enter in last equation,
since $\langle p, h_{21} \rangle = 0$.

The {\it second Lyapunov coefficient} is defined by
\begin{equation}
l_2= \frac{1}{12} \: {\rm Re} \: G_{32}, \label{defcoef2}
\end{equation}
where, from (\ref{H32m}), $G_{32}=\langle p, \mathcal H_{32}
\rangle$.

The complex vector $h_{32}$ can be found solving the nonsingular
$(n+1)$-dimensional system
\begin{equation}
\left( \begin{array}{cc}
i \omega_0 I_n -A & q \\
\\ {\bar p} & 0
\end{array} \right) \left( \begin{array} {c}
h_{32}\\
\\s \end{array} \right)= \left( \begin{array}{c} \mathcal H_{32} -G_{32} q \\
\\0 \end{array} \right),\label{h32}
\end{equation}
with the condition $\langle p, h_{32} \rangle = 0$.

From the coefficients of the terms $w^4 {\bar w}$, $w^4 {\bar
w}^2$ and $w^3 {\bar w}^3$ in (\ref{ku}), one has respectively
\begin{eqnarray}\label{h41}
h_{41}=(3 i \omega_0 I_n -A)^{-1}[4B(h_{11},h_{30})+ 6
B(h_{20},h_{21}) + 4 B(q,h_{31}) \nonumber \\ + B(\bar q,h_{40}) +
12 C(q,h_{11},h_{20}) + 6 C(q,q,h_{21}) + 4 C(q,\bar q, h_{30}) \\
+ 3 C(\bar q, h_{20},h_{20}) + 4 D(q,q,q,h_{11}) + 6 D(q,q,\bar q,
h_{20}) \nonumber \\ + E(q,q,q,q,\bar q) - 6 G_{21}h_{30}]
\nonumber,
\end{eqnarray}
\begin{eqnarray}\label{h42}
h_{42}=(2 i \omega_0 I_n -A)^{-1}[8 B(h_{11},h_{31}) + 6
B(h_{20},h_{22}) + B({\bar h}_{20},h_{40}) \nonumber \\ + 6
B(h_{21},h_{21})
+ 4 B({\bar h}_{21},h_{30}) + 4 B(q,h_{32}) + 2 B(\bar q, h_{41}) \nonumber \\
+ 12 C(h_{11},h_{11},h_{20}) + 3 C(h_{20},h_{20},{\bar h}_{20})
+ 24 C(q,h_{11},h_{21}) \nonumber \\ + 12 C(q,h_{20},{\bar h}_{21})
+ 4 C(q,{\bar h}_{20},h_{30}) + 6 C(q,q,h_{22}) + 8 C(q,\bar q,h_{31}) \nonumber \\
+ 8 C(\bar q, h_{11},h_{30}) + 12 C(\bar q, h_{20},h_{21})
+ C(\bar q, \bar q, h_{40}) \\ + 12 D(q,q,h_{11},h_{11})
+ 6 D(q,q,h_{20},{\bar h}_{20}) + 4 D(q,q,q,{\bar h}_{21}) \nonumber \\
+ 12 D(q,q,\bar q,h_{21}) + 24 D(q,\bar q,h_{11},h_{20}) + 4
D(q,\bar q,\bar q,h_{30}) \nonumber \\ + 3 D(\bar q,\bar q,
h_{20},h_{20}) + E(q,q,q,q,{\bar h}_{20}) + 8 E(q,q,q,\bar q,
h_{11}) \nonumber \\ + 6 E(q,q,\bar q,\bar q,h_{20}) +
K(q,q,q,q,\bar q,\bar q) \nonumber \\ - 4 (G_{32}h_{20} + 3
G_{21}h_{31} + {\bar G}_{21}h_{31})] \nonumber,
\end{eqnarray}
\begin{eqnarray}\label{h33}
h_{33}= - A^{-1}[9 B(h_{11},h_{22}) + 3 B(h_{20},{\bar h}_{31}) +
3 B({\bar h}_{20},h_{31}) + 9 B(h_{21},{\bar h}_{21}) \nonumber \\
+ B({\bar h}_{30},h_{30}) + 3 B(q,{\bar h}_{32}) + 3 B(\bar
q,h_{32}) + 6 C(h_{11},h_{11},h_{11}) \nonumber \\ + 9
C(h_{11},{\bar h}_{20},h_{20}) + 18 C(q,h_{11},{\bar h}_{21}) + 3
C(q,h_{20},{\bar h}_{30}) \nonumber \\ + 9 C(q,{\bar
h}_{20},h_{21}) + 3 C(q,q,{\bar h}_{31}) + 9 C(q,\bar q,h_{22}) +
18 C(\bar q,h_{11},h_{21}) \nonumber \\ + 9 C(\bar q,h_{20},{\bar
h}_{21}) + 3 C(\bar q, {\bar h}_{20},h_{30}) + 3 C(\bar q,\bar q,
h_{31}) + 9 D(q,q,{\bar h}_{20},h_{11}) \\ + D(q,q,q,{\bar
h}_{30}) + 9 D(q,q,\bar q,{\bar h}_{21}) + 18 D(q,\bar q,
h_{11},h_{11}) \nonumber \\ + 9 D(q,\bar q, {\bar h}_{20}, h_{20})
+ 9 D(q,\bar q,\bar q, h_{21}) + 9 D(\bar q,\bar q, h_{11},
h_{20}) \nonumber \\ + 3 E(q,q,q,\bar q,{\bar h}_{20}) + 9
E(q,q,\bar q,\bar q, h_{11}) + 3 E(q,\bar q,\bar q,\bar q, h_{20}) \nonumber \\
+ K(q,q,q,\bar q,\bar q, \bar q) - 3 (G_{32}+{\bar G}_{32}) h_{11}
- 9 (G_{21}+{\bar G}_{21}) h_{22}] \nonumber.
\end{eqnarray}

Defining $\mathcal H_{43} $ as
\begin{eqnarray}\label{H43m}
\mathcal H_{43} = 12 B(h_{11},h_{32}) +
6 B(h_{20},{\bar h}_{32}) + 3 B({\bar h}_{20},h_{41}) \nonumber \\
+ 18 B(h_{21},h_{22})
 + 12 B({\bar h}_{21},h_{31}) + 4 B(h_{30},{\bar
h}_{31}) + B({\bar h}_{30},h_{40}) \nonumber \\ + 4 B(q,h_{33}) +
3 B(\bar q, h_{42}) + 36 C(h_{11},h_{11},h_{21}) + 36
C(h_{11},h_{20},{\bar h}_{21}) \nonumber \\ + 12 C(h_{11},{\bar
h}_{20},h_{30}) + 3 C(h_{20},h_{20},{\bar h}_{30}) + 18
C(h_{20},{\bar h}_{20},h_{21}) \nonumber \\ + 36
C(q,h_{11},h_{22}) + 12 C(q,h_{20},{\bar h}_{31}) + 12 C(q,{\bar
h}_{20},h_{31}) \nonumber \\ + 36 C(q,h_{21},{\bar h}_{21}) + 4
C(q,h_{30},{\bar h}_{30}) + 6 C(q,q,{\bar h}_{32}) \nonumber \\ +
12 C(q,\bar q,h_{32}) + 24 C(\bar q,h_{11},h_{31}) + 18 C(\bar
q,h_{20},h_{22}) \nonumber \\ + 3 C(\bar q,{\bar h}_{20},h_{40}) +
18 C(\bar q,h_{21},h_{21}) + 12 C(\bar q,{\bar h}_{21},h_{30})
\nonumber \\ + 3 C(\bar q, \bar q, h_{41}) + 24
D(q,h_{11},h_{11},h_{11}) + 36 D(q,h_{11},h_{20},{\bar h}_{20})
\nonumber \\ + 36 D(q,q,h_{11},{\bar h}_{21}) + 6
D(q,q,h_{20},{\bar h}_{30}) + 18 D(q,q,{\bar h}_{20},h_{21})
\nonumber \\ + 4 D(q,q,q,{\bar h}_{31}) + 18 D(q,q,\bar q,h_{22})
+ 72 D(q,\bar q, h_{11},h_{21}) \\ + 36 D(q,\bar q,h_{20},{\bar
h}_{21}) + 12 D(q,\bar q,{\bar h}_{20},h_{30}) + 12 D(q,\bar
q,\bar q, h_{31}) \nonumber \\ + 36 D(\bar q,h_{11},h_{11},h_{20})
+ 9 D(\bar q,h_{20},h_{20},{\bar h}_{20}) + 12 D(\bar q,\bar
q,h_{11},h_{30}) \nonumber \\ + 18 D(\bar q,\bar q,h_{20},h_{21})
+ D(\bar q,\bar q, \bar q,h_{40}) + 12 E(q,q,q,h_{11},{\bar
h}_{20}) \nonumber \\ + E(q,q,q,q,{\bar h}_{30}) + 12 E(q,q,q,\bar
q,{\bar h}_{21}) + 36 E(q,q,\bar q,h_{11},h_{11}) \nonumber \\ +
18 E(q,q,\bar q,h_{20},{\bar h}_{20}) + 18 E(q,q,\bar q,\bar q,
h_{21}) + 36 E(q,\bar q,\bar q,h_{11},h_{20}) \nonumber \\ + 4
E(q,\bar q,\bar q,\bar q, h_{30})  + 3 E(\bar q, \bar q, \bar q,
h_{20},h_{20}) + 3 K(q,q,q,q,\bar q,{\bar h}_{20}) \nonumber \\ +
12 K(q,q,q,\bar q,\bar q, h_{11}) + 6 K(q,q,\bar q,\bar q,\bar
q,h_{20}) + L(q,q,q,q,\bar q,\bar q,\bar q) \nonumber \\ - 6
(2G_{32}h_{21} + {\bar G}_{32}h_{21} + 3 G_{21} h_{32} + 2 {\bar
G}_{21}h_{32})
\nonumber,
\end{eqnarray}
and from the coefficients of the terms $w^4 {\bar w}^3$, one has a
singular system for $h_{43}$
\begin{eqnarray}\label{h43m}
(i \omega_0 I_n -A)h_{43}= \mathcal H_{43}- G_{43} q
\end{eqnarray}
which has solution if and only if
\begin{eqnarray}\label{h43}
\langle p, \mathcal H_{43}- G_{43} q \rangle =0 ,
\end{eqnarray}
where the terms $- 6 (2G_{32}h_{21} + {\bar G}_{32}h_{21} + 3
G_{21} h_{32} + 2 {\bar G}_{21}h_{32})$ appearing in the last line
of equation (\ref{H43m}) actually do not enter in the last
equation, since $\langle p, h_{21} \rangle = 0$ and $\langle p,
h_{32} \rangle = 0$.

The {\it third Lyapunov coefficient} is defined by
\begin{equation}
l_3= \frac{1}{144} \: {\rm Re} \: G_{43}, \label{defcoef3}
\end{equation}
where, from (\ref{h43}), $ G_{43} = \langle p, \mathcal H_{43}
\rangle$.

The expressions for the vectors $h_{50}, h_{60}, h_{51}, h_{70},
h_{61}, h_{52} $ have been omitted since they are not important
here.

\begin{remark}\label{conceitual}
Other equivalent definitions and algorithmic  procedures to write
the expressions for the Lyapunov coefficients $l_j , j= 1,2,3,$
for two dimensional systems can be found in Andronov et al.
\cite{al} and Gasull et al. \cite{gt}, among others. These
procedures apply also to the three dimensional systems of this
work, if  properly restricted to the center manifold. The authors
found, however, that the  method  outlined above, due to Kuznetsov
\cite{kuznet, kuznet2}, requiring  no explicit formal evaluation
of the center manifold, is better adapted to the needs of this
work.
\end{remark}

A {\it Hopf point} $({\bf x_0}, {\bf \mu_0})$ is an equilibrium
point of (\ref{diffequat}) where the Jacobian matrix $A = f_{\bf
x}({\bf x_0}, {\bf \mu_0})$ has a pair of purely imaginary
eigenvalues $\lambda_{2,3} = \pm i \omega_0$, $\omega_0 > 0$, and
admits  no other critical eigenvalues ---i.e. located on the
imaginary axis. At a Hopf point a two dimensional center manifold
is well-defined, it is invariant under the flow generated by
(\ref{diffequat}) and can be continued with arbitrary high class
of differentiability to nearby parameter values. In fact, what is
well defined is the $\infty$-jet ---or infinite Taylor series---
of the center manifold, as well as that of its continuation, any
two of them having contact in the arbitrary high  order of their
differentiability class.

A Hopf point is called {\it transversal} if the parameter
dependent complex eigenvalues cross the imaginary axis with
non-zero derivative. In a neighborhood of a transversal Hopf point
---H1 point, for concision--- with $l_1 \neq 0$ the dynamic
behavior of the system (\ref{diffequat}), reduced to the family of
parameter-dependent continuations of the center manifold, is
orbitally topologically equivalent to the following complex normal
form
\[
w' = (\eta + i \omega) w + l_1 w |w|^2 ,
\]
$w \in \mathbb C $, $\eta$, $\omega$ and $l_1$ are real functions
having  derivatives of arbitrary  high order, which are
continuations  of $0$, $\omega_0$ and the first Lyapunov
coefficient at the H1 point. See  \cite{kuznet}. As $l_1 < 0$
($l_1 > 0$) one family of stable (unstable) periodic orbits can be
found on this family of manifolds, shrinking  to an equilibrium
point at the H1 point.

A {\it Hopf point of codimension 2} is a Hopf point where $l_1$
vanishes. It is called {\it transversal} if $\eta = 0$ and $l_1 =
0$ have transversal intersections, where $\eta = \eta (\mu)$ is
the real part of the critical eigenvalues. In a neighborhood of a
transversal Hopf point of codimension 2 ---H2 point, for
concision---  with $l_2 \neq 0$ the dynamic behavior of the system
(\ref{diffequat}), reduced to the family of parameter-dependent
continuations of the center manifold, is orbitally topologically
equivalent to
\[
w' = (\eta + i \omega_0) w + \tau w |w|^2 + l_2 w |w|^4 ,
\]
where $\eta$ and $\tau$ are unfolding parameters.  See
\cite{kuznet}. The bifurcation diagrams for $l_2 \neq 0$ can be
found in \cite{kuznet}, p. 313, and in \cite{takens}.

A {\it Hopf point of codimension 3} is a Hopf point of codimension
2 where $l_2$ vanishes. A Hopf point of codimension 3 is called
{\it transversal} if $\eta = 0$, $l_1 = 0$ and $l_2 = 0$ have
transversal intersections. In a neighborhood of a transversal Hopf
point of codimension 3 ---H3 point, for concision--- with $l_3
\neq 0$ the dynamic behavior of the system (\ref{diffequat}),
reduced to the family of parameter-dependent continuations of the
center manifold, is orbitally topologically equivalent to
\[
w' = (\eta + i \omega_0) w + \tau w |w|^2 + \nu w |w|^4 + l_3 w
|w|^6 ,
\]
where $\eta$, $\tau$ and $\nu$ are unfolding parameters. The
bifurcation diagram for $l_3 \neq 0$ can be found in Takens
\cite{takens}.

\begin{teo}
Suppose that the system
\[
{\bf x}' = f({\bf x},{\bf \mu}), \: {\bf x}=(x,y,z), \: \mu =
(\beta, \alpha, \varepsilon)
\]
has the equilibrium ${\bf x} = {\bf 0}$ for $\mu = 0$ with
eigenvalues
\[
\lambda_{2,3} (\mu) = \eta (\mu) \pm i \omega(\mu),
\]
where $\omega(0) = \omega_0 > 0$. For $\mu = 0$ the following
conditions hold
\[
\eta (0) = 0, \: l_1(0) = 0, \: l_2(0) = 0,
\]
where $l_1(\mu)$ and $l_2(\mu)$ are the first and second Lyapunov
coefficients, respectively. Assume that the following genericity
conditions are satisfied
\begin{enumerate}
\item $l_3 (0) \neq 0$, where $l_3 (0)$ is the third Lyapunov
coefficient;

\item the map $\mu \to (\eta (\mu), l_1 (\mu), l_2 (\mu))$ is
regular at $\mu = 0$.

\end{enumerate}
Then, by the introduction of a complex variable, the above system
reduced to the family of parameter-dependent continuations of the
center manifold, is orbitally topologically equivalent to
\[
w' = (\eta + i \omega_0) w + \tau w |w|^2 + \nu w |w|^4 + l_3 w
|w|^6
\]
where $\eta$, $\tau$ and $\nu$ are unfolding parameters.

\label{teoremaHopf}
\end{teo}

\begin{remark}\label{mal}
The proof of this theorem given by Takens for $C^{\infty}$
families of vector fields, using the Malgrange-Mather Preparation
Theorem \cite{gg}, is also valid in the present case of
arbitrarily high, but finite, class of differentiability, using
the appropriate extensions of the Preparation Theorem. See Bakhtin
\cite{bak} and Milman \cite {mil}, among others.
\end{remark}

\section{Hopf bifurcations}\label{hopf}

The stability of the equilibrium point $P_0$ given in
(\ref{equilibrium}) as $\varepsilon_c = \varepsilon (\beta,\alpha)
= 2 \: \alpha \: \beta ^{3/2}$ is analyzed here. According to
(\ref{taylorexp}) and the subsequent expressions (\ref{Bap}),
(\ref{Cap}), (\ref{Dap}), (\ref{Eap}), (\ref{Kap}) and
(\ref{Lap}), for $B_i$ to $L_i$, one has
\begin{equation}
A = \left( \begin{array}{ccc}
0 & 1 & 0 \\
\\- \omega_0 ^2 & - \varepsilon_c & 2 \: \beta \: \omega_0 \\
\\- \alpha \: \sqrt {\beta} \: \omega_0    & 0  & 0
\end{array} \right),
\label{partelinear}
\end{equation}
where
\begin{equation}
\omega_0 = \sqrt {\frac{1- \beta^2}{\beta}}, \label{omega0}
\end{equation}
and referring to the expressions in equations (\ref{diffequat})
and (\ref{Fhomo})
\begin{equation}
F({\bf x}) - A{\bf x} = \left( 0, F_2({\bf x}), F_3({\bf x})
\right), \label{partenaolinear}
\end{equation}
where
\begin{eqnarray*}
F_2({\bf x}) = - \frac{3}{2} \: \omega_0 \: \sqrt \beta \: x^2 +
\omega_0 \: \beta^{3/2} \: z^2 + \frac{2 (2 \beta ^2 -1)}{\sqrt
\beta} \: x \: z + \frac{4-7 \beta ^2}{6 \beta} \: x^3 \\ - 4 \:
\omega_0 \: \beta \: x^2 \: z  + (2 \beta ^2 -1) \: x \: z^2 +
\frac{5}{8} \sqrt \beta \omega_0 x^4  + \frac{4}{3} \beta^{-1/2}
(1-2\beta^2) x^3 z \\ - 2 \beta^{3/2} \omega_0 x^2 z^2 + \frac{31
\beta^2 -16}{120 \beta}  x^5 + \frac{4}{3} \beta \omega_0 x^4 z +
 \frac{2-4 \beta^2}{3} x^3 z^2 \\ - \frac{7}{80} \sqrt
\beta \omega_0 x^6 + \frac{4}{15} \beta^{-1/2} (2 \beta^2 -1) x^5
z + \frac{2}{3} \beta^{3/2} \omega_0 x^4 z^2 \\ +  \frac{64-127
\beta^2}{5040 \beta} x^7 - \frac{8}{45} \beta \omega_0 x^6 z +
 \frac{4 \beta^2 - 2}{15} x^5 z^2 + O(||{\bf
x}||^8),
\end{eqnarray*}
\begin{eqnarray*}
F_3({\bf x}) = - \frac{1}{2} \: \alpha \: \beta \: x^2 +
\frac{1}{6} \: \alpha \: \sqrt \beta \: \omega_0 x^3 +
\frac{1}{24} \: \alpha \: \beta \: x^4 - \frac{1}{120} \: \alpha
\: \sqrt \beta \: \omega_0 x^5 \\ - \frac{1}{720} \: \alpha \:
\beta \: x^6 + \frac{1}{5040} \: \alpha \: \sqrt \beta \: \omega_0
x^7 + O(||{\bf x}||^8).
\end{eqnarray*}
\noindent From equations (\ref{taylorexp}), (\ref{Bap}),
(\ref{Cap}), (\ref{Dap}), (\ref{Eap}), (\ref{Kap}), (\ref{Lap})
and (\ref{partenaolinear}) one has
\begin{equation}
B({\bf x},{\bf y}) = \left( 0, B_2({\bf x},{\bf y}), -\alpha \:
\beta \: x_1 \: y_1 \right), \label{B1}
\end{equation}
\noindent where
\begin{eqnarray*}
B_2({\bf x},{\bf y})= -3 \: \omega_0 \: \sqrt \beta \: x_1 \: y_1
+ 2 \: \omega_0 \: \beta^{3/2} \: x_3 \: y_3 + \frac{2 \: (2 \:
\beta ^2 -1)}{\sqrt \beta} \: ( x_1 \: y_3 + x_3 \: y_1),
\end{eqnarray*}
\begin{equation}
C({\bf x},{\bf y}, {\bf z}) = \left( 0, C_2({\bf x},{\bf y},{\bf
z}),\alpha \: \sqrt \beta \: \omega_0 \: x_1 \: y_1 \: z_1
 \right),
\label{C1}
\end{equation}
\noindent where
\begin{eqnarray*}
C_2({\bf x},{\bf y},{\bf z})= \frac{4-7 \beta^2}{\beta} \: x_1 \:
y_1 \: z_1 - 8 \omega_0 \: \beta \: (x_1 \; y_1 \; z_3 + x_1 \;
y_3 \; z_1 + x_3 \; y_1 \; z_1 ) \\ + 2 \: (2 \beta^2 -1) \: (x_1
\; y_3 \; z_3 + x_3 \; y_1 \; z_3 + x_3 \; y_3 \; z_1),
\end{eqnarray*}
\begin{equation}
D({\bf x},{\bf y}, {\bf z},{\bf u}) = \left( 0, D_2({\bf x},{\bf
y},{\bf z},{\bf u}), \alpha \beta x_1 y_1 z_1 u_1 \right),
\label{D1}
\end{equation}
\noindent where
\begin{eqnarray*}
D_2({\bf x},{\bf y},{\bf z},{\bf u}) = 15 \omega_0 \beta^{1/2} x_1
y_1 z_1 u_1 + 8 \left( \frac{1-2 \beta^2}{\beta^{1/2}} \right) \:
\Big( x_1 y_1 z_1 u_3 + x_1 y_1 z_3 u_1 \\ + x_1 y_3 z_1 u_1 + x_3
y_1 z_1 u_1 \Big) - 8 \omega_0 \beta^{3/2} \Big( x_1 y_1 z_3 u_3 +
x_1 y_3 z_1 u_3 \\ + x_1 y_3 z_3 u_1 + x_3 y_3 z_1 u_1 + x_3 y_1
z_1 u_3 + x_ 3 y_1 z_3 u_1 \Big),
\end{eqnarray*}
\begin{equation}
E({\bf x},{\bf y}, {\bf z},{\bf u},{\bf v}) = \left( 0, E_2({\bf
x},{\bf y},{\bf z},{\bf u},{\bf v}), - \alpha \omega_0 \beta^{1/2}
x_1 y_1 z_1 u_1 v_1 \right), \label{E1}
\end{equation}
\noindent where
\begin{eqnarray*}
E_2({\bf x},{\bf y},{\bf z},{\bf u},{\bf v}) = \frac {31 \beta^2 -
16}{\beta} \; x_1 y_1 z_1 u_1 v_1 + 32 \omega_0 \beta \: \Big( x_1
y_1 z_1 u_1 v_3 \\ + x_1 y_1 z_1 u_3 v_1 + x_1 y_1 z_3 u_1 v_1 +
x_1 y_3 z_1 u_1 v_1 + x_3 y_1 z_1 u_1 v_1 \Big) \\ + 8 (1 - 2
\beta^2) \Big( x_1 y_1 z_1 u_3 v_3 + x_1 y_1 z_3 u_1 v_3 + x_1 y_1
z_3 u_3 v_1 \\ + x_1 y_3 z_1 u_1 v_3 + x_1 y_3 z_1 u_3 v_1 + x_1
y_3 z_3 u_1 v_1 + x_3 y_3 z_1 u_1 v_1 \\ + x_3 y_1 z_3 u_1 v_1 +
x_3 y_1 z_1 u_3 v_1 + x_3 y_1 z_1 u_1 v_3 \Big),
\end{eqnarray*}
\begin{equation}
K({\bf x},{\bf y}, {\bf z},{\bf u},{\bf v},{\bf s}) = \left( 0,
K_2({\bf x},{\bf y},{\bf z},{\bf u},{\bf v},{\bf s}),  -\alpha
\beta x_1 y_1 z_1 u_1 v_1 s_1 \right), \label{K1}
\end{equation}
\noindent where
\begin{eqnarray*}
K_2({\bf x},{\bf y},{\bf z},{\bf u},{\bf v},{\bf s})= -63 \omega_0
\beta^{1/2} x_1 y_1 z_1 u_1 v_1 s_1 + 32 \left( \frac{2 \beta^2 -
1}{\beta^{1/2}} \right) \Big( x_1 y_1 z_1 u_1 v_1 s_3 \\+ x_1 y_1
z_1 u_1 v_3 s_1 + x_1 y_1 z_1 u_3 v_1 s_1 + x_1 y_1 z_3 u_1 v_1
s_1 + x_1 y_3 z_1 u_1 v_1 s_1 \\ + x_3 y_1 z_1 u_1 v_1 s_1 \Big) +
32 \omega_0 \beta^{3/2} \Big( x_1 y_1 z_1 u_1 v_3 s_3 + x_1 y_1
z_1 u_3 v_1 s_3 + x_1 y_1 z_3 u_1 v_1 s_3 \\ + x_1 y_3 z_1 u_1 v_1
s_3 + x_3 y_1 z_1 u_1 v_1 s_3 + x_1 y_1 z_1 u_3 v_3 s_1 + x_1 y_1
z_3 u_1 v_3 s_1 \\ + x_1 y_3 z_1 u_1 v_3 s_1 + x_3 y_1 z_1
u_1 v_3 s_1 + x_1 y_1 z_3 u_3 v_1 s_1 + x_1 y_3 z_1 u_3 v_1 s_1 \\
+ x_3 y_1 z_1 u_3 v_1 s_1 + x_1 y_3 z_3 u_1 v_1 s_1 + x_3 y_1 z_3
u_1 v_1 s_1 + x_3 y_3 z_1 u_1 v_1 s_1 \Big),
\end{eqnarray*}
\begin{equation}
L({\bf x},{\bf y}, {\bf z},{\bf u},{\bf v},{\bf s},{\bf t}) =
\left( 0, L_2({\bf x},{\bf y}, {\bf z},{\bf u},{\bf v},{\bf
s},{\bf t}),  \omega_0 \beta^{1/2} x_1 y_1 z_1 u_1 v_1 s_1 t_1
\right), \label{L1}
\end{equation}
\noindent where
\begin{eqnarray*}
L_2({\bf x},{\bf y},{\bf z},{\bf u},{\bf v},{\bf s},{\bf t})= (64
\omega_0 ^2 - 63 \beta) x_1 y_1 z_1 u_1 v_1 s_1 t_1 - 128 \omega_0
\beta \Big( x_1 y_1 z_1 u_1 v_1 s_1 t_3 \\ + x_1 y_1 z_1 u_1 v_1
s_3 t_1 + x_1 y_1 z_1 u_1 v_3 s_1 t_1 + x_1 y_1 z_1 u_3 v_1 s_1
t_1 + x_1 y_1 z_3 u_1 v_1 s_1 t_1 + \\ x_1 y_3 z_1 u_1 v_1 s_1 t_1
+ x_3 y_1 z_1 u_1 v_1 s_1 t_1 \Big) + 32 \beta (\beta - \omega_0
^2) \Big( x_1 y_1 z_1 u_1 v_1 s_3 t_3 + \\ x_1 y_1 z_1 u_1 v_3 s_1
t_3 + x_1 y_1 z_1 u_3 v_1 s_1 t_3 +  x_1 y_1 z_3 u_1 v_1 s_1 t_3 +
x_1 y_3 z_1 u_1 v_1 s_1 t_3 + \\ x_3 y_1 z_1 u_1 v_1 s_1 t_3 + x_1
y_1 z_1 u_1 v_3 s_3 t_1 + x_1 y_1 z_1 u_3 v_1 s_3 t_1 + x_1 y_1
z_3 u_1 v_1 s_3 t_1 + \\ x_1 y_3 z_1 u_1 v_1 s_3 t_1 + x_3 y_1 z_1
u_1 v_1 s_3 t_1 + x_1 y_1 z_1 u_3 v_3 s_1 t_1 + x_1 y_1 z_3 u_1
v_3 s_1 t_1 + \\ x_1 y_3 z_1 u_1 v_3 s_1 t_1 + x_3 y_1 z_1 u_1 v_3
s_1 t_1 + x_1 y_1 z_3 u_3 v_1 s_1 t_1 + x_1 y_3 z_1 u_3 v_1 s_1
t_1 + \\ x_3 y_1 z_1 u_3 v_1 s_1 t_1 + x_1 y_3 z_3 u_1 v_1 s_1 t_1
+ x_3 y_1 z_3 u_1 v_1 s_1 t_1 + x_3 y_3 z_1 u_1 v_1 s_1 t_1 \Big).
\end{eqnarray*}

The eigenvalues of $A$ (equation (\ref{partelinear})) are
\begin{equation}
\lambda_1= -\varepsilon_c = -2 \alpha \beta^{3/2}, \: \: \lambda_2
= i \: \omega_0, \: \: \lambda_3 = -i \: \omega_0.
\label{autovalores}
\end{equation}
and from (\ref{normalization}) one has
\begin{equation}
q = \left( -i, \omega_0, \frac{\varepsilon_c}{2 \beta} \right)
\label{q}
\end{equation}
\noindent and
\begin{equation}
p = \left( -\frac{i}{2}, \frac{\omega_0 - i \varepsilon_c}{2 (
\omega_0 ^2 + \varepsilon_c ^2 )}, \frac{\beta (\varepsilon_c + i
\omega_0)}{ \omega_0 ^2 + \varepsilon_c ^2 } \right).
\label{p}
\end{equation}

\begin{figure}[!h]

\centerline{
\includegraphics[width=10cm]{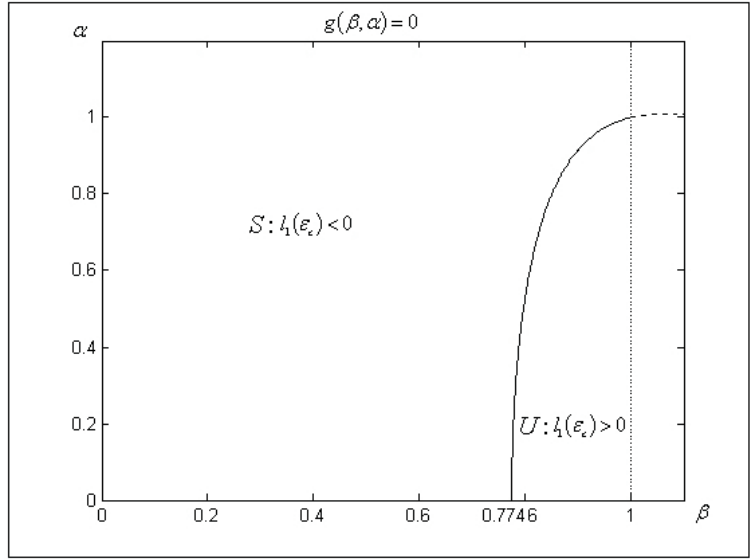}}

\caption{{\small Signs of the first Lyapunov coefficient}.}

\label{curvag=0}

\end{figure}

\begin{teo}
Consider the three-parameter family of differential equations
(\ref{wattdem}). The first Lyapunov coefficient is given by
\begin{equation}
l_1 (\beta,\alpha,\varepsilon_c) = - \frac{1}{2} \: \left(
\frac{\alpha \beta^{3/2} (1 - \beta ^2) \left( 3 + (\alpha^2 -5)
\beta ^2 + \alpha ^4 \beta ^6 \right)}{ \left( 1 - \beta^2 +
\alpha^2 \beta^4 \right) \left( 1 - \beta^2 + 4 \alpha^2 \beta^4
\right)} \right). \label{coeficiente}
\end{equation}
If
\begin{equation}
g(\beta, \alpha) =  3 + (\alpha^2 -5) \beta ^2 + \alpha ^4 \beta
^6 \label{hopfcond}
\end{equation}
\noindent is different from zero then the three-parameter family
of differential equations (\ref{wattdem}) has a transversal Hopf
point at $P_0$ for $\varepsilon_c = \varepsilon (\beta,\alpha) = 2
\: \alpha \: \beta ^{3/2}$.

If $(\beta, \alpha, \varepsilon_c) \in S \cup U$ (see Fig.
\ref{curvag=0}) then the three-parameter family of differential
equations (\ref{wattdem}) has a H1 point at $P_0$. If $(\beta,
\alpha, \varepsilon_c) \in S$ then the H1 point at $P_0$ is
asymptotically stable and for each $\varepsilon < \varepsilon_c$,
but close to $\varepsilon_c$, there exists a stable periodic orbit
near the unstable equilibrium point $P_0$. If $(\beta, \alpha,
\varepsilon_c) \in U$ then the H1 point at $P_0$ is unstable and
for each $\varepsilon > \varepsilon_c$, but close to
$\varepsilon_c$, there exists an unstable periodic orbit near the
asymptotically stable equilibrium point $P_0$.

\label{teohopf}
\end{teo}

This theorem  summarizes Proposition 3.2 and Theorems 3.5, 3.6 and
3.7 established in \cite{smb}. Equation (\ref{hopfcond}) gives a
simple expression for the sign of the first Lyapunov coefficient
(\ref{coeficiente}). Its graph is illustrated in Fig.
\ref{curvag=0}, where the signs of the first Lyapunov coefficient
are also represented. The curve $l_1 = 0$ divides the surface of
critical parameters into two connected components denoted by $S$
and $U$ where $l_1 < 0$ and $l_1 > 0$ respectively.

We have the following theorem.
\begin{teo}
Consider the three-parameter family of differential equations
(\ref{wattdem}) restricted to $\varepsilon =\varepsilon_c $. The
second Lyapunov coefficient is given by
\begin{equation}
l_2 (\beta,\alpha,\varepsilon_c) =  \frac{\alpha \: \beta^{3/2} \:
h(\beta,\alpha,\varepsilon_c)}{ 36(1 - \beta^2 + \alpha^2
\beta^4)^3(9 - 9 \beta^2 + 4 \alpha^2 \beta^4) (1- \beta^2 + 4
\alpha^2 \beta^4)^3}, \label{coeficiente2}
\end{equation}
where
\begin{eqnarray*}
h(\beta,\alpha,\varepsilon_c)= -162 - 54(-9 + 37 \alpha^2) \beta^2
- 9 (-126 + 61 \alpha^2 + 60 \alpha^4) \beta^4 - 18(405 \\ - 3212
\alpha^2 + 1128 \alpha^4) \beta^6 + (13770 - 210843 \alpha^2 +
113612 \alpha^4 - 5533 \alpha^6) \beta^8 \\ - 6 (2133 - 57687
\alpha^2 + 38218 \alpha^4 + 5186 \alpha^6) \beta^{10} \\ + (5994 -
301275 \alpha^2 + 215340 \alpha^4 + 284264 \alpha^6 - 16022
\alpha^8) \beta^{12}
\\ + 2(-567 + 67878 \alpha^2 - 45196 \alpha^4 - 379430 \alpha^6 +
9347 \alpha^8) \beta^{14} \\ + \alpha^2( - 25029 + 9540 \alpha^2 +
990831 \alpha^4 + 155856 \alpha^6 - 21205 \alpha^8) \beta^{16} \\
+ 4 \alpha^4 (513 - \alpha^2 (163340 + 120616 \alpha^2 - 16768
\alpha^4)) \beta^{18} \\ - 2 \alpha^6 (-86887 - 258835 \alpha^2 +
30173 \alpha^4 + 7208 \alpha^6) \beta^{20} + 2 \alpha^8 (-96867 \\
- 8956 \alpha^2 + 23208 \alpha^4) \beta^{22}
+ \alpha^{10}(33671 - 58288 \alpha^2 - 4880 \alpha^4) \beta^{24} \\
+ 16 \alpha^{12} (1603 + 718 \alpha^2) \beta^{26} - 16 \alpha^{14}
(453 + 40 \alpha^2) \beta^{28} + 640 \alpha^{16} \beta^{30}.
\end{eqnarray*}

\label{lemacoef2}
\end{teo}

\noindent{\bf Proof.} Define the following functions
\[
T_1 = {\rm Re} \langle p, E(q,q,q,\bar q, \bar q) \rangle, \: T_2
= {\rm Re} \langle p, D(q,q,q,{\bar h}_{20}) \rangle,
\]
\[
T_3 = {\rm Re} \langle p, D(q, \bar q,\bar q, h_{20}) \rangle, \:
T_4 = {\rm Re} \langle p, D(q,q,\bar q,h_{11}) \rangle,
\]
\[
T_5 = {\rm Re} \langle p, C(\bar q, \bar q, h_{30}) \rangle, \:
T_6 = {\rm Re} \langle p, C(q,q,{\bar h}_{21}) \rangle, \: T_7 =
{\rm Re} \langle p, C(q,\bar q, h_{21}) \rangle,
\]
\[
T_8 = {\rm Re} \langle p, C(q, {\bar h}_{20}, h_{20}) \rangle, \:
T_9 = {\rm Re} \langle p, C(q,h_{11},h_{11}) \rangle,
\]
\[
T_{10} = {\rm Re} \langle p, C(\bar q, h_{20}, h_{11}) \rangle, \:
T_{11} = {\rm Re} \langle p, B(\bar q, h_{31}) \rangle, \: T_{12}
= {\rm Re} \langle p, B(q,h_{22}) \rangle,
\]
\[
T_{13} = {\rm Re} \langle p, B({\bar h}_{20},h_{30}) \rangle, \:
T_{14} = {\rm Re} \langle p, B({\bar h}_{21},h_{20}) \rangle, \:
T_{15} = {\rm Re} \langle p, B(h_{11},h_{21}) \rangle.
\]
From (\ref{defcoef2}) one has
\begin{eqnarray*}
{\rm Re} \: G_{32} = T_1 + T_2 + 3 T_3 + 6 T_4 + T_5 + 3 T_6 + 6
T_7 + 3 T_8 + 6 T_9 + 6 T_{10} \\ + 2 T_{11} + 3 T_{12} + T_{13} +
3 T_{14} + 6 T_{15}.
\end{eqnarray*}

The theorem follows by expanding the expressions in definition of
the second Lyapunov coefficient (\ref{defcoef2}). It relies on
extensive calculation involving the vector $q$ (\ref{q}), the
vector $p$ (\ref{p}), the functions $B$, $C$, $D$ and $E$, listed
equations (\ref{B1}), (\ref{C1}), (\ref{D1}) and (\ref{E1}),
respectively, the long complex vectors $h_{11}$, $h_{20}$,
$h_{30}$, $h_{21}$, $h_{31}$ and $h_{22}$, and the above functions
$T_1$ to $T_{15}$.

The calculations in this proof, corroborated by Computer Algebra,
have been posted in \cite{mello}. Here, the complex vectors
$h_{11}$, $h_{20}$, $h_{30}$, $h_{21}$, $h_{31}$ and $h_{22}$ have
particularly  long expressions. They have listed in MATHEMATICA 5
files.
\begin{flushright}
$\blacksquare$
\end{flushright}

\begin{figure}[!h]

\centerline{
\includegraphics[width=10cm]{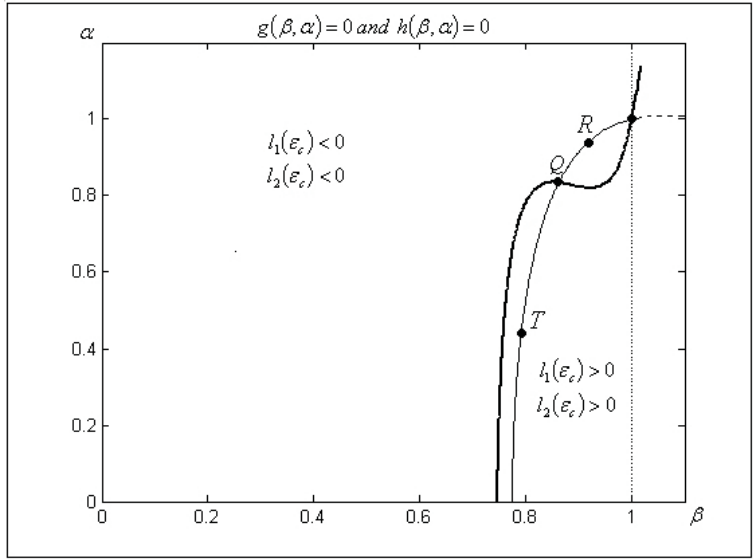}}

\caption{{\small Signs of the first and second Lyapunov
coefficients}.}

\label{curvagh}

\end{figure}

\begin{teo}
For the system (\ref{wattdem}) there is unique point $Q = (\beta,
\alpha, \varepsilon_c)$, with coordinates
\[
\beta = 0.86828033997971281542...,\, \alpha =
0.85050048430685017856...
\]
and
\[
\varepsilon_c = 1.37624106484659953171...
\]
where the curves $l_1 = 0$ and $l_2 = 0 $ on the critical surface
intersect and there do it transversally.

\label{crucial2}
\end{teo}

\noindent{\bf Computer assisted Proof.} The point $Q$ is the
intersection of the curves $l_1 = 0 $ and $l_2 = 0$ on the Hopf
critical surface. It is defined and obtained by the solution of
the equations
\[
g(\beta, \alpha) = 0,
\]
given in (\ref{hopfcond}), and
\begin{equation}
h(\beta, \alpha) = h(\beta, \alpha, \varepsilon_c) = 0,
\label{hopfcond2}
\end{equation}
where $h(\beta, \alpha, \varepsilon_c)$ is given by
(\ref{coeficiente2}). The existence and uniqueness of $Q$ with the
above coordinates has been established numerically with the
software MATHEMATICA 5.

\begin{figure}[!h]

\centerline{
\includegraphics[width=10cm]{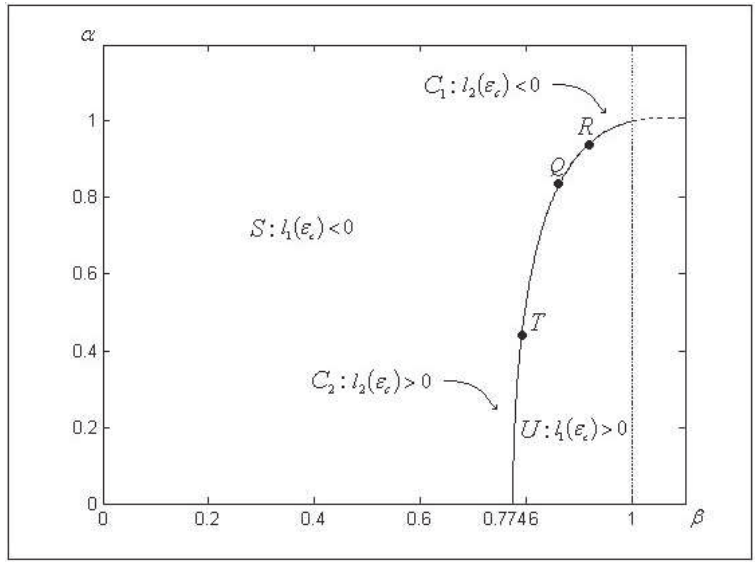}}

\caption{{\small Signs of $l_2$ on the curve $l_1 = 0$}.}

\label{pointsTR}

\end{figure}

Figure \ref{curvagh} presents a geometric synthesis interpreting
the long calculations involved in this proof. The sign of
$h(\beta, \alpha)$ gives the sign of the second Lyapunov
coefficient (\ref{coeficiente2}). The graph of $h(\beta,\alpha) =
0$, where the signs of the first and second Lyapunov coefficients
are also illustrated. As follows, $l_2 < 0$ on the open arc of the
curve $l_1 = 0$, denoted  by $C_1$. On this arc a typical
reference point $R$ is depicted. Also $l_2 > 0$ on the open arc of
the curve $l_1= 0$, denoted by $C_2$. This arc  contains the
typical reference point, denoted by $T$. See also Fig.
\ref{pointsTR}.

\begin{figure}[!h]
\centerline{ {\includegraphics[width=12cm]{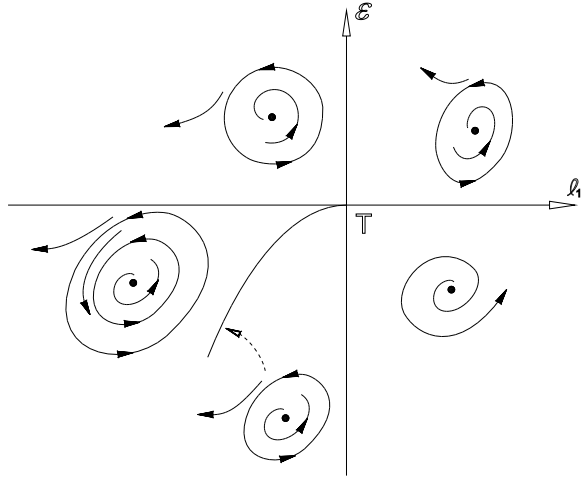}}}

\caption{{\small Bifurcation diagram of the system (\ref{wattdem})
at point $T$}.}

\label{pointT}

\end{figure}

\begin{figure}[!h]
\centerline{ {\includegraphics[width=12cm]{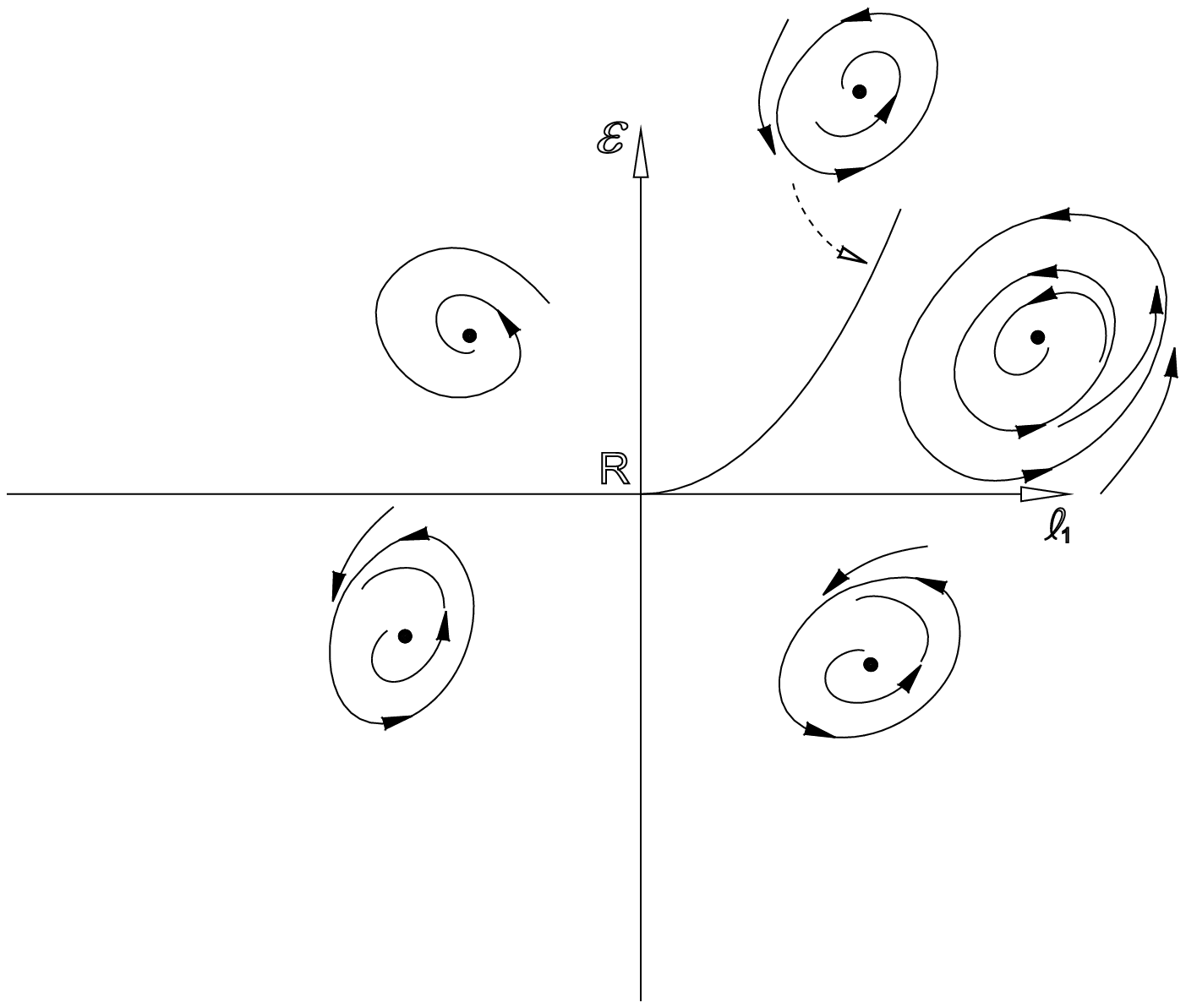}} }

\caption{{\small Bifurcation diagram of the system (\ref{wattdem})
at point $R$}.}

\label{pointR}

\end{figure}

The bifurcation diagrams of the system (\ref{wattdem}) at the
points $T$ and $R$ are illustrated in Fig. \ref{pointT} and
\ref{pointR}, as a consequence of \cite{kuznet} and \cite{takens}.

The main steps of the calculations  that provide the numerical
evidence for this theorem  have been posted in \cite{mello}.
\begin{flushright}
$\blacksquare$
\end{flushright}

\begin{teo}
If $(\beta, \alpha, \varepsilon_c) \in C_1 \cup C_2$ then the
three-parameter family of differential equations (\ref{wattdem})
has a transversal Hopf point of codimension 2 at $P_0$. If
$(\beta, \alpha, \varepsilon_c) \in C_2$ then the H2 point at
$P_0$ is unstable and the bifurcation diagram is drawn  in Fig.
\ref{pointT}. If $(\beta, \alpha, \varepsilon_c) \in C_1$ then the
H2 point at $P_0$ is asymptotically stable and the bifurcation
diagram is illustrated in Fig. \ref{pointR}.

\label{teohopf2}
\end{teo}

This theorem is a synthesis of the discussion in the last part in
the proof of Theorem \ref{crucial2}.

\begin{figure}[!h]

\centerline{
\includegraphics[width=10cm]{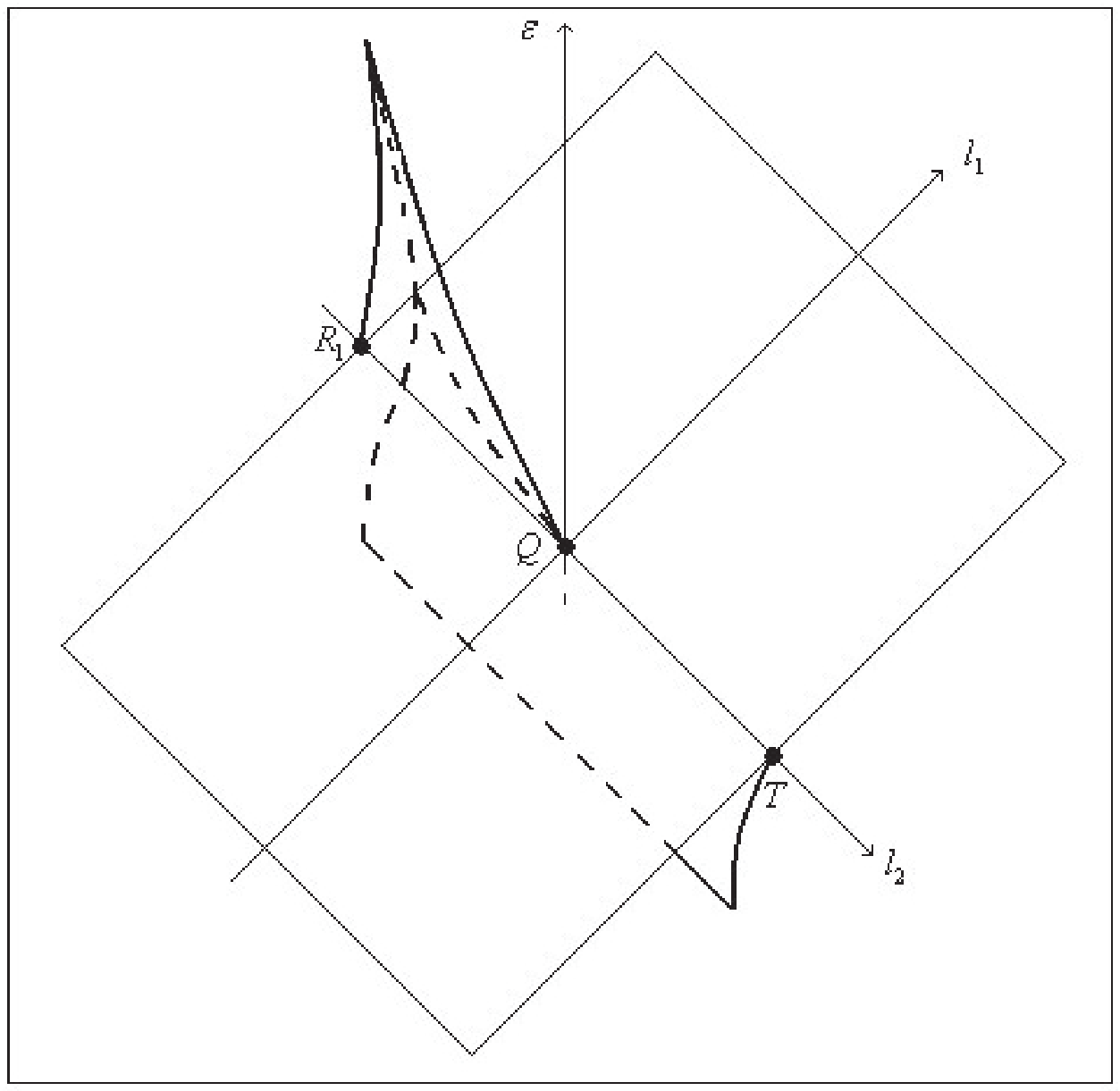}}

\caption{{\small Bifurcation diagram of the system (\ref{wattdem})
at point $Q$}.}

\label{pointQ}

\end{figure}

\begin{figure}[!h]

\centerline{
\includegraphics[width=13cm]{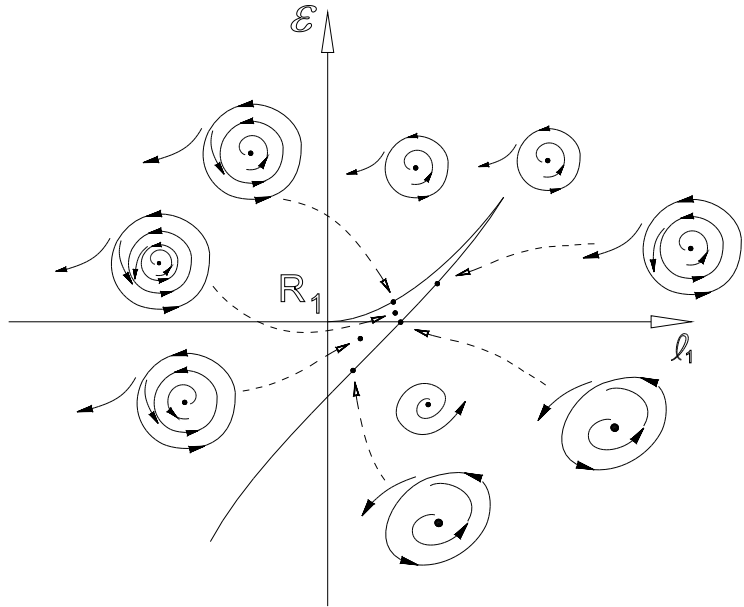}}

\vspace{-1cm}

\caption{{\small Bifurcation diagram of the system (\ref{wattdem})
at point $R_1$}.}

\label{pointR1}

\end{figure}

\begin{teo}
For the parameter values at the point $Q$ determined in Theorem
\ref{crucial2}, the three-parameter family of differential
equations (\ref{wattdem}) has a tranversal Hopf point of
codimension 3 at $P_0$ which is asymptotically unstable since
$l_3(Q) > 0$. The bifurcation diagram of system (\ref{wattdem}) at
the point $Q$ is illustrated in Figs. \ref{pointQ} and
\ref{pointR1}.

\label{teohopf3}
\end{teo}

\noindent{\bf Computer Assisted  Proof.} For the  point $Q$ take
five decimal round-off coordinates $\beta = 0.86828$, $\alpha =
0.85050$ and $\varepsilon_c = 1.37624$. For these values of the
parameters one has
\[
p = \left( -i/2, 0.12224 - 0.31601 i, 0.54878 + 0.21228 i \right),
\]
\[
q = \left( -i, 0.53237, 0.79250 \right),
\]
\[
h_{11}= \left( -1.75030, 0, 0.48792 \right),
\]
\[
h_{20}= \left( -2.24198 - 0.11191 i, 0.11916 -
 2.38715 i, 0.04434 - 1.58196 i \right),
\]
\[
h_{30}= \left( -2.68329 + 5.27951 i, -8.43202 -
 4.28554 i, -4.24045 - 0.86409 i \right),
\]
\begin{equation}
G_{21}= - 2.90053 i, \label{G21l3}
\end{equation}
\[
h_{21}= \left( 1.20918 + 0.65492 i, -3.24920 + 0.64374 i, 1.26042
+ 1.11353 i \right),
\]
\[
h_{40}= \left( 9.27690 + 25.24802 i, -53.76550 + 19.75510 i,
-9.11345 + 11.36572 i \right),
\]
\[
h_{31}= \left( -25.72175 - 5.12199 i, 4.47976 - 7.87822 i, 6.22842
- 15.97687 i \right),
\]
\[
h_{22} = \left( -15.72589, 0, 10.92671 \right),
\]
\begin{equation}
G_{32} = -34.93331 i, \label{G32l3}
\end{equation}
\[
h_{32} = \left( 27.17768 + 53.16361 i, -57.53733 + 3.94677 i,
52.73722 + 27.89259 i \right),
\]
\[
h_{41} = \left( -35.5370 + 180.2333 i, -195.9736 - 10.0589 i,
-125.3480 - 33.7428 i \right),
\]
\[
h_{42} = \left( -778.4924 - 466.4510 i, 362.1612 + 81.2385 i,
390.2364 - 503.3807 i \right),
\]
\[
h_{33} = \left( -536.09324, 0, 835.33555 \right),
\]
\begin{equation}
G_{43} = 56.23254 - 2424.27069 i .\label{G43l3}
\end{equation}
From (\ref{defcoef}), (\ref{defcoef2}), (\ref{defcoef3}),
(\ref{G21l3}), (\ref{G32l3}) and (\ref{G43l3}) one has
\[
l_1 (Q) = 0, \:\: l_2 (Q) = 0, \:\: l_3 (Q) = \frac{1}{144} \:
{\rm Re} \: G_{43} = 0.39050.
\]
The calculations above have also been corroborated with 20
decimals round-off precision performed using the software
MATHEMATICA 5 \cite{math}. See \cite{mello}.

The gradients of the functions $l_1 $, given in
(\ref{coeficiente}), and $l_2 $, given in (\ref{coeficiente2}), at
the point $Q$ are, respectively
\[
(0.80095, -0.31847), (-0.38861, -0.85118).
\]
The transversality condition at $Q$ is equivalent to the
non-vanishing of the determinant of the matrix whose columns are
the above gradient vectors, which is evaluated gives $-0.80552$.
The transversality condition being satisfied, the bifurcation
diagrams in Figs. \ref{pointQ} and \ref{pointR1}, follow from the
work of Takens \cite{takens}, taking into consideration the
orientation and signs established in Theorems \ref{crucial2} and
\ref{teohopf2}.
\begin{flushright}
$\blacksquare$
\end{flushright}

\section{Concluding comments}\label{conclusion}

The historical relevance of the Watt governor study as well as its
importance for present day theoretical and technological aspects
of Automatic Control has been widely discussed by Denny
\cite{denny} and others. See also \cite{mac, smb}.

This paper starts reviewing  the  stability analysis due to
Maxwell and Vyshnegradskii, which accounts for the
characterization, in the space of parameters, of the structural as
well as Lyapunov stability of the equilibrium of the Watt
Centrifugal Governor System, WGS. It continues with recounting the
extension of the analysis  to the first order, codimension one
stable points, happening on the complement of a curve in the
critical surface where the eigenvalue criterium of Lyapunov holds,
as studied in \cite{has1}, \cite{humadi} and by the authors
\cite{smb}, based on the calculation of the first Lyapunov
coefficient. Here the bifurcation analysis at the equilibrium
point of the WGS is pushed forward  to the calculation of the
second and third Lyapunov coefficients which make possible the
determination of the Lyapunov as well as higher order structural
stability at the equilibrium point.  See also \cite{kuznet,
kuznet2}, \cite{gt} and \cite{al} .

The calculations of these coefficients, being extensive, rely  on
Computer Algebra and Numerical evaluations carried out with the
software MATHEMATICA 5 \cite{math}. In the site \cite{mello} have
been posted the main steps of the calculations in the form of
notebooks for MATHEMATICA 5.

With the analytic and numeric data provided in the analysis
performed here, the bifurcation diagrams are established along the
points of the curve where the first Lyapunov coefficient vanishes.
Picture \ref{pointQ} and \ref{pointR1} provide  a qualitative
synthesis of the  dynamical conclusions achieved  here at the
parameter values where the WGS achieves  most complex equilibrium
point. A reformulation of these conclusions follow:

There is a ``solid tongue"  where two stable regimes coexist: one
is an equilibrium and the other is  a small amplitude periodic
orbit, i.e. an oscillation.

For parameters inside the ``tongue", this conclusion suggests, a
{\it hysteresis} explanation for the phenomenon of ``hunting"
observed in the performance of WGS in an early stage of the
research on its stability conditions. Which attractor represents
the actual state of the system will depend on the path along which
the parameters evolve to reach their actual values of the
parameters under consideration. See Denny \cite{denny} for
historical comments, where he refers to the term ``hunting" to
mean an oscillation around an equilibrium going near but not
reaching it.

Finally, we would like to stress that although this work
ultimately focuses the specific three dimensional, three parameter
system of differential equations given by (\ref{wattde}), the
method of analysis and calculations explained in Section
\ref{codim} can be adapted to the study of other systems with
three or more phase variables and depending on three or more
parameters.

\vspace{0.2cm}

\noindent {\bf Acknowledgement}: The first and second authors
developed this work under the project CNPq Grant 473824/04-3. The
first author is fellow of CNPq and takes part in the project CNPq
PADCT 620029/2004-8. This work was finished while he visited Brown
University, supported by FAPESP, Grant 05/56740-6.

The authors are grateful to C. Chicone and R. de la Llave for
helpful comments.

\end{document}